\documentclass{amsart}

\usepackage{amssymb,amsfonts,amsxtra}
\usepackage[all]{xypic}
\xyoption{dvips}

\theoremstyle{definition}
\newtheorem{ntn}{Notation}

\theoremstyle{plain}
\newtheorem{lem}[ntn]{Lemma}
\newtheorem{prp}[ntn]{Proposition}
\newtheorem{thm}[ntn]{Theorem}
\newtheorem{cor}[ntn]{Corollary}

\theoremstyle{remark}

\newtheorem{exa}[ntn]{Example}

\newcommand{\m}{\mathfrak{m}}
\newcommand{\xymat}{\SelectTips{cm}{}\xymatrix}
\newcommand{\C}{\mathbb{C}}
\newcommand{\D}{\mathcal{D}}
\newcommand{\I}{\mathcal{I}}
\newcommand{\N}{\mathbb{N}}
\renewcommand{\O}{\mathcal{O}}
\newcommand{\V}{\mathcal{V}}
\newcommand{\Z}{\mathbb{Z}}

\DeclareMathOperator{\depth}{depth}
\DeclareMathOperator{\gr}{gr}
\DeclareMathOperator{\grade}{grade}
\DeclareMathOperator{\img}{im}
\DeclareMathOperator{\rk}{rk}
\DeclareMathOperator{\syz}{syz}
\DeclareMathOperator{\Ann}{Ann}
\DeclareMathOperator{\Der}{Der}
\DeclareMathOperator{\DR}{DR}
\DeclareMathOperator{\GL}{GL}
\DeclareMathOperator{\LL}{L}
\DeclareMathOperator{\RR}{R}
\DeclareMathOperator{\Rees}{Rees}
\DeclareMathOperator{\Sing}{Sing}
\DeclareMathOperator{\SC}{Sp}
\DeclareMathOperator{\Sym}{Sym}
\DeclareMathOperator{\T}{T}

\begin{document}

\title[A criterion for the logarithmic differential operators]{A criterion for the logarithmic differential operators to be generated by vector fields}
\author{Mathias Schulze}
\address{Mathias Schulze\\
Department of Mathematics\\
Oklahoma State University\\
401 MSCS\\
Stillwater, OK 74078\\
USA}
\email{mschulze@math.okstate.edu}

\begin{abstract}
We study divisors in a complex manifold in view of the property that the algebra of logarithmic differential operators along the divisor is generated by logarithmic vector fields.
We give 
\begin{itemize}
\item a sufficient criterion for the property,
\item a simple proof of F.J. Calder{\'o}n--Moreno's theorem that free divisors have the property,
\item a proof that divisors in dimension $3$ with only isolated quasi--homogeneous singularities have the property,
\item an example of a non--free divisor with non--isolated singularity having the property,
\item an example of a divisor not having the property, and
\item an algorithm to compute the V--filtration along a divisor up to a given order.
\end{itemize}
\end{abstract}



\thanks{The author is grateful to M.~Granger for many valuable discussions and comments and to F.J.~Castro--Jim\'enez, L.~Narv\'aez--Macarro, and J.M.~Ucha--Enr\'iquez for explaining their results and ideas.}

\maketitle
\tableofcontents

\section{Logarithmic comparison theorem for free divisors}

Let $X$ be a complex manifold of dimension $n\ge2$, $\O$ the ring holomorphic functions on $X$, and $\Omega^\bullet$ the complex of holomorphic differential forms.
Grothendieck's Comparison Theorem states that the De~Rham system $\O$ is regular \cite[Thm.\,2.3.4]{Meb89}.
This is equivalent to the fact that, for any divisor $D\subset X$, the natural morphism
\[
\xymat{
\Omega^\bullet(*D)=\DR(\O(*D))\ar[r]&\RR i_*i^{-1}\DR(\O)=\RR i_*\C_U,
}
\]
where $i$ is the inclusion $U=X\backslash D\subset X$, is a quasi--isomorphism.
Let  $\Omega^\bullet(\log D)\subset\Omega^\bullet(*D)$ be the subcomplex of logarithmic differential forms along $D$ \cite[Def.\,1.2]{Sai80}.
The above statement raises the question whether the inclusion $\Omega^\bullet(\log D)\subset\Omega^\bullet(*D)$ is also a quasi--isomorphism.
That is: Can one compute the cohomology of the complement of $D$ by logarithmic differential forms along $D$?
This turns out to be a property of $D$ called the logarithmic comparison theorem or simply LCT.
It is an open problem to characterize the divisors for which LCT holds.

Let $\Theta=\Der_\C(\O)$ be the $\O$--module of holomorphic vector fields on $X$ and $\Der(\log D)\subset\Theta$ be the $\O$--submodule of logarithmic differential operators along $D$ \cite[Def.\,1.4]{Sai80}.
A divisor $D$ is called free if $\Der(\log D)$ is a locally free $\O$--module.
Let $\D$ be the $\O$--algebra of differential operators on $X$ with holomorphic coefficients and let $F$ be the increasing filtration on $\D$ by the order of differential operators.
Let $\V^D$ be the V--filtration along $D$ on $\D$ as defined in Section \ref{20} such that  $\V^D_0=\D(\log D)$ is the $\O$--algebra of logarithmic differential operators along $D$.
F.J.~Calder{\'o}n--Moreno \cite[Thm.\,1]{Cal99} proves that, for a free divisor $D$,  $\V^D_0$ is generated by vector fields, that is $\V^D_0=\O[\Der(\log D)]$.
Let $S_D$ be the decreasing filtration on $\D$ which is locally defined by $S_D^k=f^k\cdot\D$ where $f\in\O$ such that $D=(f)$.
By Corollary \ref{24}, the induced filtration $S_D$ on $\V^D_0$ defined by $S_D^k\V^D_0=\V^D_0\cap( f^k\cdot\D)$ reflects the embeddings $\V^D_k\subset\D$.
If $\V^D_0=\O[\Der(\log D)]$ then $(\V^D_0,S_D)$ is a filtered $(\V^D_0,\V^D)$--module.
 
F.J.~Calder{\'o}n--Moreno and L.~Narv{\'a}ez--Macarro \cite[Cor.\,4.2]{CN05} prove that LCT holds for a free divisor $D$ if and only if the complex
\[
\D\overset\LL\otimes_{\D(\log D)}\O(D)=\D\otimes_{\D(\log D)}\SC^\bullet_{\D(\log D)}(\O(D))
\]
is concentrated in degree $0$ and the natural multiplication morphism 
\[
\xymat{
\D\otimes_{\D(\log D)}\O(D)\ar[r]^-{\epsilon_D}&\O(*D)
}
\]
is an isomorphism.
Injectivity of $\epsilon_D$ is locally equivalent to $\Ann_\D(\frac{1}{f})$ being generated by operators of order $1$ where $f\in\O$ such that $D=(f)$.
For any divisor $D$, T.~Torrelli proves that the latter condition already implies surjectivity of $\epsilon_D$ \cite[Prop.\,1.3]{Tor04} and conjectures that it is even equivalent to LCT \cite[Conj.\,1.11]{Tor04}.
A problem to verify this conjecture for a free divisor $D$ consists in $\D\otimes_{\D(\log D)}\SC^\bullet_{\D(\log D)}(\O(D))$ not being $F$--strict  in general \cite[Rem.\,4.2.4]{Cal99}.
So grading by $F$ does not reduce the problem to a commutative one.
But both properties of $D$ in question can be characterized in terms of $S_D$--strictness:
On the one hand, exactness of $\D\otimes_{\D(\log D)}\SC^\bullet_{\D(\log D)}(\O(D))$ in degree $k$ is equivalent to $S_D$--strictness of the differential of $\SC^\bullet_{\D(\log D)}(\O(D))$ from degree $k-1$ to degree $k$.
On the other hand, injectivity of $\epsilon_D$ is equivalent to $S_D$--strictness of the last differential of $\D\otimes_{\D(\log D)}\SC^\bullet_{\D(\log D)}(\O(D))$.

A solution of the LCT problem seems to require a deeper understanding of the V--filtration in general.
There are many questions:
\begin{itemize}
\item What are properties of the $\V^D_k$?
\item When is $\V^D_0$ generated by vector fields?
\item When is $\V^D_0$ locally finitely generated?
\item What are properties of the embeddings $\V^D_k\subset\D$?
\end{itemize}
We shall approach the first two questions in this article.

\section{V--filtration along subvarieties and divisors}\label{20}

Let $Y\subset X$ be a subvariety in $X$ and let $\I\subset\O$ be its ideal.
The V--filtration $\V^Y$ along $Y$ is the increasing filtration on $\D$ defined by
\[
\V^Y_k=\{P\in\D\mid\forall l\in\Z:P(\I^l)\subset\I^{l-k}\}
\]
for all $k\in\Z$.
We shall omit the index $Y$ if it is clear from the context.
Clearly $\V_k\cdot\V_l\subset\V_{k+l}$ for all $k,l\in\Z$.
Hence $\V_0$ is an $\O$--algebra and $\V_k$ is an $\V_0$--module for all $k\in\Z$. 

\begin{exa}\label{2}
Let $x_1,\dots,x_m,y_1,\dots,y_n$ be coordinates on $X=\C^{m+n}$.
\begin{enumerate}
\item For the submanifold $Y=\{y=0\}$,
\[
\V^Y_k=\bigl\{P=\sum_{j_1-i_1+\cdots+j_n-i_n\le k}P_{i,j}(x,\partial_x)y_1^{i_1}\partial_{y_1}^{j_1}\cdots y_n^{i_n}\partial_{y_n}^{j_n}\in\D\bigr\}.
\]
\item For the normal crossing divisor $D=(y_1\cdots y_n)$,
\[
\V^D_k=\bigl\{P=\sum_{j_1-i_1,\dots,j_n-i_n\le k}P_{i,j}(x,\partial_x)y_1^{i_1}\partial_{y_1}^{j_1}\cdots y_n^{i_n}\partial_{y_n}^{j_n}\in\D\bigr\}.
\]
\end{enumerate}
\end{exa}

Denote the complement of the singularities of $Y$ by
\[
\xymat{
U_Y=X\backslash\Sing(Y)\ar[r]^-{i_Y}&X.
}
\]
We shall omit the index $Y$ if it is clear from the context.
The V--filtration along a divisor has a special property.

\begin{prp}\label{1}
Let $D\subset X$ be a divisor.
Then $\V^D=(i_D)_*i_D^{-1}\V^D$.
\end{prp}

\begin{proof}
We may assume that $D=(f)$ for some $f\in\O$ by the local nature of the statement.
Since $\V_k\subset\D$ and $\D$ is a locally free $\O$--module, 
\[
i_*i^{-1}\V_k\subset i_*i^{-1}\D=\D.
\]
Since $\O\cdot f^{l-k}$ is a free $\O$--module, $P\in i_*i^{-1}\V_k$ implies
\[
P(g\cdot f^l)\in i_*i^{-1}(\O\cdot f^{l-k})=\O\cdot f^{l-k}
\]
for all $g\in\O$ and $l\in\Z$ and hence $P\in\V_k$.
\end{proof}

\begin{cor}\label{24}
Let $D=(f)\subset X$ with $f\in\O$ be a divisor.
Then
\[
\V_k=
\begin{cases}
f^{-k}\V_0, & k\le0,\\
f^{-k}(\V_0\cap f^k\D), & k\ge1.
\end{cases}
\]
\end{cor}

\begin{proof}
The equalities in question hold on $U_D$ by Example \ref{2} (2) and hence on $X$ by Proposition \ref{1}.
\end{proof}

Denote the symbol map for $F$ by
\[
\xymat{
\D\ar@{->>}[r]^-\sigma&\gr^F\D
}.
\]
The decomposition $F_1\D=\O\oplus\Theta$ defines the $\O$--module $\Der(\log Y)\subset\Theta$ of logarithmic vector fields along $Y$ by
\[
F_1\V_0=\O\oplus\Der(\log Y).
\]
This definition simplifies to
\[
\Der(\log Y)=\{\theta\in\Theta\mid \theta(\I)\subset\I\}
\]
by the Leibniz rule and implies involutivity of $\Der(\log Y)$, that is
\[
[\Der(\log Y),\Der(\log Y)]\subset\Der(\log Y).
\]

\begin{exa}\label{3}
Let $x_1,\dots,x_m,y_1,\dots,y_n$ be coordinates on $X=\C^{m+n}$.
\begin{enumerate}
\item For the submanifold $Y=\{y=0\}$,
\[
\Der(\log D)=\O\langle\partial_{x_1},\dots,\partial_{x_m}\rangle+\O\langle y_i\partial_{y_j}\mid1\le i,j\le n\rangle.
\]
\item For the normal crossing divisor $D=(y_1\cdots y_n)$, 
\[
\Der(\log D)=\O\langle\partial_{x_1},\dots,\partial_{x_m},y_1\partial_{y_1},\dots,y_n\partial_{y_n}\rangle.
\]
\end{enumerate}
\end{exa}

Let $\O[\Der(\log Y)]\subset\D$ be the image of the tensor algebra
\[
\xymat{
\T_\C\Der(\log Y)\ar[r]^-{\gamma_Y}&\D.
}
\]
Then at least $\O[\Der(\log Y)]\subset\V_0^Y$.

\begin{cor}\label{7}
Let $D\subset X$ be a divisor.
Then $\V_0^D=\O[\Der(\log D)]$ if and only if $\O[\Der(\log D)]=(i_D)_*i_D^{-1}\O[\Der(\log D)]$.
\end{cor}

\begin{proof}
By Examples \ref{2} and \ref{3}, $\V_0^D=\O[\Der(\log D)]$ on $U_D$.
Hence the claim follows from Proposition \ref{1}.
\end{proof}

A divisor $D\subset X$ is called free if $\Der(\log D)$ is a locally free $\O$--module.
By K.~Saito \cite[Cor.\,1.7]{Sai80}, $\Der(\log D)$ is reflexive and hence all divisors in dimension $n=2$ are free.
By Example \ref{3} (2), normal crossing divisors are free.
In particular, any divisor $D$ is free on $U_D$.
 
F.J.~Calder\'on--Moreno \cite[Thm.\,1]{Cal99} proves that $\V^D_0=\O[\Der(\log D)]$ for a free divisor.
We give a simple proof of this result.

\begin{cor}
Let $D\subset X$ be a free divisor.
Then $\V^D_0=\O[\Der(\log D)]$.
\end{cor}

\begin{proof}
By Lemma \ref{23} and grading by $F$, $\O[\Der(\log D)]$ is a locally free $\O$--module and hence Corollary \ref{7} applies.
\end{proof}

\begin{lem}\label{23}
Let $R$ be a domain and let $P_1,\dots,P_n\in R\cdot T_1\oplus\dots\oplus R\cdot T_n=R^n$ be $R$--linearly independent.
Then $R[P_1,\dots,P_n]\subset R[T_1,\dots,T_n]$ is a polynomial ring.
\end{lem}

\begin{proof}
Write $P_i=\sum_j p_{i,j}T_j$ with $p_{i,j}\in R$.
Then by assumption $p=\det(p_{i,j})\ne0$ and hence $R_p[P_1,\dots,P_n]$ is a polynomial ring.
Since $R$ is a domain, $\xymat{R\ar[r]&R_p}$ is injective and hence $R[P_1,\dots,P_n]$ is a polynomial ring.
\end{proof}

In general it is not clear if, or under which conditions, $\V_0^Y$ is a locally finite $\O$--algebra.
Even to compute $F_k\V_0^Y$ is a problem since the definition involves infinitely many conditions.
The following result allows one to compute $F_k\V_0^D$ algorithmically.

\begin{prp}\label{5}
Let $x_1,\dots,x_n$ be coordinates on $X=\C^n$.
Let $D=(f)\subset X$ with $f\in\O$ be a divisor.
Then, for $P\in F_d\D$, $P\in\V^D_k$ if and only if 
\begin{equation}\label{4}
\forall\alpha\in\N^n,l\in\N:|\alpha|+l\le d\Rightarrow P(x^\alpha f^l)\in\O\cdot f^{l-k}.
\end{equation}
\end{prp}

\begin{proof}
Let $0\ne P\in F_d\D$ and assume that condition (\ref{4}) holds.
For $l\in\N$, the vector space $\C[x_1,\dots,x_n]_{\le d-l}$ is invariant under $x\mapsto Ax+a$ for $a\in\C^n$ and $A\in\GL_n(\C)$.
Hence, at a smooth point $y$ of $D$, condition (\ref{4}) holds for coordinates $x_1,\dots,x_n$ at $y$ such that $\partial_{x_n}(f)(y)\ne0$.
Then $y_1,\dots,y_{n-1},t=x_1,\dots,x_{n-1},f$ are coordinates at $y$ such that 
\[
\forall\beta\in\N^{n-1},l\in\N:|\beta|+l\le d\Rightarrow P_y(y^\beta t^l)\in\O_y\cdot t^{l-k}.
\]
Write $P_y=\sum_{|\beta|+l\le d}p_{\beta,l}\partial_y^\beta\partial_t^l$ with $p_{\beta,l}\in\O_y$ and choose $\gamma\in\N^{n-1}$ and $m\in\N$ such that $|\gamma|+m$ is minimal with $p_{\gamma,m}\ne0$.
Then
\[
\gamma!m!p_{\gamma,m}=P(y^\gamma t^m)\in\O_y\cdot t^{m-k}
\]
and hence $p_{\gamma,m}\partial_y^\gamma\partial_t^m\in\V_{k,y}$ by Example \ref{2} (2).
By increasing induction on $|\gamma|+m$, this implies $P_y\in\V_{k,y}$ for all $y\in U_D$ and hence $P\in\V_k$ by Proposition \ref{1}.
\end{proof}

\begin{exa}\label{9}
Let $x,y,z$ be coordinates on $\C^3$ and
\[
f=xyz(x+y+z)(x+2y+3z).
\]
Then $D=(f)\subset\C^3$ is a central generic hyperplane arrangement.
Let
\[
Q=(x+y+z)(x+2y+3z)(3zy^2\partial_y^2+(x+4y-3z)yz\partial_y\partial_z-4yz^2\partial_z^2).
\]
Then $Q\in F_2\V^D_0$ by a {\sc Singular} \cite{GPS05} computation using Proposition \ref{5}.
We shall see in Example \ref{16} that $Q\notin F_2\O[\Der(\log D)]$.
\end{exa}

There is another special property of the V--filtration along a divisor.

\begin{prp}
Let $D\subset X$ be a divisor.
Then $\depth_x(\V^D_k)\ge2$ for all $x\in X$ and $k\in\Z$.
\end{prp}

\begin{proof}
Let $x\in X$ and $D_x=(f)$ with $f\in\O_x$.
Since $\O_x$ is torsion free and $\depth(\O_x)\ge2$, there is an $\O_x$--sequence $a_1,a_2\in\m_x$ such that $a_1$ is different from all irreducible factors of $f$.
Let $P\in\V_{k,x}$ with $a_2\cdot P\in a_1\cdot\V_{k,x}\subset a_1\cdot\D_x$.
Then $P\in a_1\cdot\D_x$ since $\D_x$ is a free $\O_x$--module.
But $P(g\cdot f^l)\in\O_x\cdot f^{l-k}$ implies $(a_1^{-1}\cdot P)(g\cdot f^l)\in\O_x\cdot f^{l-k}$ by the choice of $a_1$ for all $g\in\O$ and $l\in\Z$ and hence $P\in a_1\cdot\V_{0,x}$.
Then $a_1,a_2\in\m_x$ is a $\V_{k,x}$--sequence and hence $\depth_x(\V_k)\ge2$.
\end{proof}

\section{Symmetric algebra of logarithmic vector fields}
 
The condition in Corollary \ref{7} is difficult to verify in general.
Therefore we focus on a case in which it still holds after grading by $F$.
There is a commutative diagram of graded algebras
\[
\xymat@C=-20pt{
&\T_\C\Der(\log Y)\ar@{->>}[ld]_-{\Sigma}\ar[rd]^-{\gr\gamma_Y}\\
\Sym_\O\Der(\log Y)\ar[rr]^-{\alpha_Y}\ar@{->>}[d]_-{\pi_Y}&&\gr^F\O[\Der(\log Y)]\\
\Rees_\O\Der(\log Y)\ar[rr]^\cong&&\O[\sigma(\Der(\log Y)]\ar@{^(->}[u]
}.
\]

\begin{lem}\label{12}
If $\alpha_Y$ is an isomorphism then
\[
\Sym_\O\Der(\log Y)=(i_Y)_*i_Y^{-1}\Sym_\O\Der(\log Y)
\]
implies $\O[\Der(\log Y)]=(i_Y)_*i_Y^{-1}\O[\Der(\log Y)]$.
\end{lem}

\begin{proof}
There is a commutative diagram
\[
\xymat@C=-40pt{
\Sym_\O\Der(\log Y)\ar[rr]\ar[d]^-{\alpha_Y}&&i_*i^{-1}\Sym_\O\Der(\log Y)\ar[d]^-{i_*i^{-1}\alpha_Y}\\
\gr^F\O[\Der(\log Y)]\ar@{^(->}[rr]\ar@{^(->}[rd]&&i_*i^{-1}\gr^F\O[\Der(\log Y)]\\
&\gr^Fi_*i^{-1}\O[\Der(\log Y)]\ar@{^(->}[ru]
}.
\]
Then the claim follows by induction on $\deg(P)$ for $P\in i_*i^{-1}\O[\Der(\log Y)]$.
\end{proof}

\begin{lem}\label{10}
$\alpha_Y$ is an isomorphism if and only if $\pi_Y$ is injective.
\end{lem}

\begin{proof}
Assume that $\pi_Y$ is injective.
An element of $\gr^F\O[\Der(\log Y)]$ is of the form $\sigma(\gamma_Y(P))$ where $P\in\T_\C\Der(\log Y)$.
Write $P=P_0\oplus\cdots\oplus P_d$ where $d=\deg(P)$.
If $\sigma(\gamma_Y(P))\notin\img\alpha_Y$ then $(\gr\gamma_Y)(P_d)=(\gr\gamma_Y)(P)=0$ and hence $P_d\in\ker\Sigma$ by injectivity of $\pi_Y$.
By definition of $\Sym_\O$, this implies that $P_d$ is in the two--sided ideal generated by the relations $\xi\otimes\eta-\eta\otimes\xi$ and $\xi\otimes(a\eta)-(a\xi)\otimes\eta$ where $\xi,\eta\in\Der(\log Y)$ and $a\in\O$. 
But
\[
\gamma_Y(\xi\otimes\eta-\eta\otimes\xi)=\xi\eta-\eta\xi=[\xi,\eta]\in\Der(\log Y)
\]
by involutivity of $\Der(\log Y)$ and 
\[
\gamma_Y(\xi\otimes(a\eta)-(a\xi)\otimes\eta)=\xi a\eta-a\xi\eta=[\xi,a]\eta=\xi(a)\eta\in\Der(\log Y).
\]
This means that
\begin{align*}
\xi\otimes\eta-\eta\otimes\xi&\equiv[\xi,\eta]\mod\ker\gamma_Y,\\
\deg(\xi\otimes\eta-\eta\otimes\xi)&>\deg([\xi,\eta]),\\
\xi\otimes(a\eta)-(a\xi)\otimes\eta&\equiv\xi(a)\eta\mod\ker\gamma_Y,\\
\deg(\xi\otimes(a\eta)-(a\xi)\otimes\eta)&>\deg(\xi(a)\eta).
\end{align*}
Hence $\gamma_Y(P)=\gamma_Y(P')$ and $\deg(P)<\deg(P')$ for some $P'\in\T_\C\Der(\log Y)$.
Then the claim follows by induction on $d=\deg(P)$.
\end{proof}

\begin{exa}\label{16}
Let $D$ and $Q$ be as in Example \ref{9}.
Then a {\sc Singular} \cite{GPS05} computation shows that $\pi_D$ is injective and that
\[
\sigma(Q)\notin\alpha_D\bigl(\Sym_\O^2\Der(\log D)\bigr).
\]
By Lemma \ref{10}, this implies $Q\notin F_2\O[\Der(\log D)]$ and hence, by Example \ref{9}, $\O[\Der(\log D)]\subsetneq\V^D_0$.
\end{exa}

By the following general statement, injectivity of $\pi_Y$ is equivalent to $\O$--torsion freeness of $\Sym_\O\Der(\log Y)$.

\begin{lem}\label{11}
Let $R$ be a domain and $M$ a finitely presented torsion free $R$--module.
Then the following are equivalent:
\begin{enumerate}
\item $\Sym_RM$ is $R$--torsion free.
\item $\Sym_RM$ is a domain.
\item $\xymat{\Sym_RM\ar@{->>}[r]^-{\pi_M}&\Rees_RM}$ is injective.
\end{enumerate}
\end{lem}

\begin{proof}
Assume that $\Sym_RM$ is $R$--torsion free.
Let $K=Q(R)$ be the fraction field of $R$.
Then $M\otimes_RK\cong K^d$ where $d=\rk(M)$.
By choosing a basis of $K^d$ and clearing denominators, one can embed $M\subset R^d$.
Then
\[
\Sym_R(M)\otimes_RK\cong\Sym_{R\otimes_RK}(M\otimes_RK)\cong\Sym_K(K^d)
\]
is a domain and hence $\Sym_RM$ is a domain since $R$ is a domain.
Applying $\Sym_R$ to the inclusion $M\subset R^d$ yields
\[
\xymat@C=-10pt{
\Sym_RM\ar[rr]^-\phi\ar@{->>}[rd]_-{\pi_M}&&\Sym_R(R^d)\\
&\Rees_RM\ar@{^(->}[ru]
}
\]
and $\ker(\phi)\otimes_RK=0$ since
\[
\Sym_R(M)\otimes_RK\cong\Sym_K(K^d)\cong\Sym_R(R^d)\otimes_RK.
\]
A presentation 
\[
\xymat{
R^m\ar[r]^-{(a_{i,j})}&R^n\ar[r]&M\ar[r]&0
}
\]
of $M$ defines an isomorphism
\[
\Sym_RM\cong R[T_1,\dots,T_n]/J
\]
where $J=\langle\sum_ja_{i,j}T_j\rangle$ is a prime ideal since $\Sym_RM$ is a domain.
Since $\Sym_R(R^d)$ is a domain, $\ker\phi$ lifts to a prime ideal $Q\subset R[T_1,\dots,T_n]$.
Then $J\subset Q$, $Q\cap R=0$, and $J\otimes_RK=Q\otimes_RK$ implies $J=Q$ and hence $\ker\pi_M=\ker\phi=0$.
\end{proof}

\begin{exa}
Let $D_4\subset\C^4$ be the central generic hyperplane arrangement defined in Section \ref{18}.
Then one can compute that the coordinates are zero divisors on $\Sym_\O^2A_4$ at $0$.
By Lemmata \ref{10}, \ref{11}, and \ref{17}, this implies that $\alpha_{D_4}$ is not an isomorphism.
\end{exa}

A divisor $D\subset X$ is called Euler homogeneous if locally $\chi(f)=f$ for some $\chi\in\Der(\log D)$ and $f\in\O$ such that $D=(f)$.
In this case, $\chi$ is called an Euler vector field and
\[
\Der(\log D)\cong\O\cdot\chi\oplus\Ann_\Theta(f).
\]
If $\Der(\log D)\cong\O\cdot\chi\oplus A$ then $\Sym_\O\Der(\log D)\cong\Sym_\O(A)[\chi]$.
For an Euler homogeneous divisor $D$, this implies
\[
\Sym_\O\Der(\log D)\cong\Sym_\O(\Ann_\Theta(f))[\chi].
\]

\begin{prp}\label{13}
Let $D\subset X$ be a divisor such that $\Sym_\O\Der(\log D)$ is $\O$--torsion free.
Then $\V^D_0=\O[\Der(\log D)]$ follows from
\[
\Sym_\O\Der(\log D)=(i_D)_*i_D^{-1}\Sym_\O\Der(\log D).
\]
If $D$ is Euler homogeneous and $A=\Ann_\Theta(f)$ or $A\oplus\O\cdot\chi\cong\Der(\log D)$ then the latter is equivalent to $\Sym_\O A=(i_D)_*i_D^{-1}\Sym_\O A$.
\end{prp}

\begin{proof}
This follows from Corollary \ref{7}, Lemmata \ref{12}, \ref{10}, and \ref{11}, and the preceding remarks.
\end{proof}

\section{Depth and torsion of symmetric algebras}

Using a theorem of G.~Scheja \cite{Sch61} on extension of coherent analytic sheaves, we shall give sufficient conditions for $\V^D_0=\O[\Der(\log D)]$ in terms of the depth and torsion of the symmetric algebras in Proposition \ref{13}.

\begin{thm}\label{14}
Let $D\subset X$ be a divisor such that $\Sym_\O\Der(\log D)$ is $\O$--torsion free.
Let $Z\subset\Sing(D)$ be a closed subset such that $\V^D_0=\O[\Der(\log D)]$ on $X\backslash Z$.
Then $\V^D_0=\O[\Der(\log D)]$ on $X$ if
\[
\depth_z\bigl(\Sym^k_\O\Der(\log D)\bigr)\ge\dim_z(Z)+2
\]
for all $z\in Z$ and $k\in\N$.
In particular, this holds if $D$ is Euler homogeneous, $A=\Ann_\Theta(f)$ or $A\oplus\O\cdot\chi\cong\Der(\log D)$, and
\[
\depth_z\bigl(\Sym^k_\O A\bigr)\ge\dim_z(Z)+2
\]
for all $z\in Z$ and $k\in\N$.
\end{thm}

\begin{proof}
This follows from \cite[Satz\,I-III]{Sch61} and Proposition \ref{13}.
\end{proof}

We shall apply a criterion by C.~Huneke \cite{Hun81} for the torsion freeness of symmetric algebras. 

\begin{prp}\label{22}
Let $R$ be a Noetherian domain and let
\[
\xymat@C=20pt{
0\ar[r]&R\ar[rr]^-{(a_1,\dots,a_m)^t}&&R^m\ar[r]&M\ar[r]&0
}
\]
be a resolution of $M$.
If $\grade(I)\ge k+1$ for $I=\langle a_1,\dots,a_m\rangle$ then
\[
\depth(I,\Sym_R(M))\ge k.
\]
\end{prp}

\begin{proof}
We proceed by induction on $k$.
By \cite[Prop.\,2.1]{Hun81}, $\grade(I)\ge 2$ implies that $\Sym_R(M)$ is $R$--torsion free.
If $k\ge2$ then $\grade(I/a)\ge k$ for some $a\in I$ and hence
\[
\xymat{
0\ar[r]&R/a\ar[rr]^-{([a_1],\dots,[a_m])^t}&&(R/a)^m\ar[r]&M/a\ar[r]&0
}
\]
is a resolution of $M$.
Since $\Sym_{R/a}(M/a)\cong\Sym_R(M)/a$, the induction hypothesis applies.
\end{proof}

\begin{thm}\label{21}
Let $D\subset X$ be an Euler homogeneous divisor and $A=\Ann_\Theta(f)$ or $A\oplus\O\cdot\chi\cong\Der(\log D)$.
Let $Z\subset\Sing(D)$ be a closed subset such that $\V^D_0=\O[\Der(\log D)]$ on $X\backslash Z$.
For $z\in Z$, let
\[
\xymat{
0\ar[r]&\O_z\ar[rr]^-{(a_{z,1},\dots,a_{z,m})^t}&&\O_z^m\ar[r]&A_z\ar[r]&0
}
\]
be a resolution of $A_z$ such that
\[
\grade(\langle a_{z,1},\dots,a_{z,m}\rangle)\ge\dim_z(Z)+3.
\]
Then $\V^D_0=\O[\Der(\log D)]$ on $X$.
\end{thm}

\begin{proof}
This follows from Theorem \ref{14}, Proposition \ref{22}, and \cite[Prop.\,2.1]{Hun81}.
\end{proof}

\begin{cor}
Let $X$ be a complex manifold of dimension $3$ and let $D\subset X$ be a divisor with only isolated quasi--homogeneous singularities.
Then $\V^D_0=\O[\Der(\log D)]$.
\end{cor}

\begin{proof}
We may assume that $X\subset\C^3$ is an open neighbourhood of $0$, $D=(f)$ with $f\in\O$, and $\Sing(D)=\{0\}$.
Let $x_1,x_2,x_3$ be coordinates on $X$.
Then $\partial(f)=\partial_1(f),\partial_2(f),\partial_3(f)\in\m_0$ is an $\O_0$--sequence.
Hence the Koszul--complex
\[
\xymat@C=10pt@R=10pt{
0\ar[r]&\O_0\ar[rr]^-{\partial(f)^t}&&\O_0^3\ar[rr]\ar[rd]&&\O_0^3\ar[rr]&&\O_0\ar[rr]&&\O_0/\langle\partial(f)\rangle\ar[r]&0\\
&&&&\Ann_{\Theta_0}(f)\ar[ru]\ar[rd]\\
&&&0\ar[ru]&&0
}
\]
is exact and induces a resolution of $\Ann_{\Theta_0}(f)$.
Then the claim follows from Theorem \ref{21} with $Z=\Sing(D)$.
\end{proof}

Our criterion also applies to some cases of non--isolated singularities.

\begin{exa}
Let $D_3\subset\C^3$ be the central generic hyperplane arrangement defined in Section \ref{18}.
Then $D_3$ is not a free divisor and has a non--isolated singularity at $0$.
By Lemma \ref{19} and Proposition \ref{17}, $A_3\cong\O^3/\O\cdot(x_1,x_2,x_3)$ and $A_3\oplus\O\cdot\chi\cong\Der(\log D_3)$.
Then, by Examples \ref{2} (2) and \ref{3} (2) on $\Sing(D_3)\backslash\{0\}$ and Theorem \ref{21} for $A=A_3$ and $Z=\{0\}$, $\V^{D_3}_0=\O[\Der(\log D_3)]$.
\end{exa}

Our approach may fail in dimension $n>3$ even for quasi--homogeneous isolated singularities.

\begin{exa}
Let $x_1,x_2,x_3,x_4$ be coordinates on $\C^4$ and 
\[
f=x_1^2+x_2^2+x_3^2+x_4^2.
\]
Then $D=(f)\subset\C^4$ has a quasi--homogeneous isolated singularity at $0$.
One can compute that the coordinates are zero divisors on $\Sym_\O^2\Ann_\O(f)$ at $0$.
By Lemmata \ref{10} and \ref{11} this implies that $\alpha_{D}$ is not an isomorphism.
\end{exa}

\section{Example of generic hyperplane arrangements}\label{18}

We shall provide some background for the examples in the previous sections.
Let $x_1,\cdots,x_n$ be coordinates on $\C^n$ and 
\[
f_n=x_1\cdots x_n(x_1+\cdots+x_n).
\]
Then $D_n=(f_n)\subset\C^n$ is a central generic hyperplane arrangement.
Let $\chi=\sum_ix_i\partial_i$ be the Euler vector field, 
\[
\eta_{i,j}=x_ix_j(\partial_i-\partial_j)\in\Der(\log D_n)
\]
for $i<j$, and $A_n=\O\langle\eta_{i,j}\rangle$.
By J.~Wiens \cite[Thm.\,3.4]{Wie01},
\[
\Der(\log D_n)=\O\cdot\chi+A_n
\]
with a minimal number of generators.
Let
\[
\sigma_{i,j,k}=x_i\eta_{j,k}-x_j\eta_{i,k}+x_k\eta_{i,j}\in\syz(\eta_{i,j})
\]
for $i<j<k$ and choose a monomial ordering refining $\partial_1<\cdots<\partial_n$.

\begin{lem}\label{19}
$(\eta_{i,j})$ is a standard basis of $A_n$ and $\syz(\eta_{i,j})=\langle\sigma_{i,j,k}\rangle$.
\end{lem}

\begin{proof}
This follows from Buchberger's criterion \cite[Thm.\,1.7.3]{GP02}.
\end{proof}

\begin{prp}\label{17}
$\Der(\log D_n)=\O\cdot\chi\oplus A_n$.
\end{prp}

\begin{proof}
It suffices to verify that no syzygy of $\chi$ and the $\eta_{i,j}$ involves $\chi$.
One can obtain the syzygies from a standard basis computation \cite[Alg.\,2.5.4]{GP02}.
The first s--polynomials $x_k\chi-\eta_{k,n}$ and $x_j\eta_{i,k}-x_i\eta_{j,k}$ have a zero $\partial_n$ component.
Hence only a sequence of s--polynomials starting with $x_k\chi-\eta_{k,n}$ can contribute to syzygies involving $\chi$ and the coefficient of $\chi$ remains a monomial.
Each element in such a sequence has exactly one monomial involving $x_n$.
Since the $\partial_2,\dots,\partial_{n-1}$ are leading components of the $\eta_{i,j}$, the sequence terminates with a non--zero element $a_k\partial_1\equiv x^{\alpha_k}\chi\mod A_n$.
By the same reason, $\O\cdot\partial_1\oplus A_n$ is a direct sum and hence $x^{\alpha_j}a_k=x^{\alpha_k}a_j$.
This implies that the coefficient of $\chi$ in any syzygy is zero.
\end{proof}

\bibliographystyle{amsalpha}
\bibliography{cldo}

\end{document}